\documentclass[11pt,a4paper,reqno]{amsart}
\usepackage{amsfonts}
\usepackage{graphicx}   
\usepackage{graphics}
\usepackage{epsfig}
\usepackage{amssymb}
\usepackage{amsmath}
\usepackage{anysize}
\marginsize{2.5cm}{2.5cm}{2.5cm}{2.5cm}

\newtheorem{theorem}{Theorem}[section]

\newtheorem{lemma}{Lemma}[section]

\newtheorem{corollary}{Corollary}[section]
\newtheorem{definition}{Definition}[section]

\newcounter{Ex} 
\newcounter{Example-ENM}
\newcounter{Example-A} 
\newcounter{Example-B} 

\def\({\left(}
\def\){\right)}

\newlength\savedwidth



\begin{document}

\setcounter{page}{0}
\title [N.A. Bokayev et al.: Cones of monotone functions ]
{Cones of monotone functions generated by a generalized fractional maximal function}

\author[TWMS J. Pure Appl. Math., V.15, N.1, 2024] {N.A. Bokayev$^{1}$, A.Gogatishvili$^{2}$, A.N. Abek$^{1,*}$}
\thanks{$^{1}$L.N. Gumilyov Eurasian National University, K.Munaitpasov str., 23, 010008 Astana, Kazakhstan\\
 \indent \,\,\, e-mail: bokayev2011@yandex.ru, azhar.abekova@gmail.com$^{*}$ \\
 \indent $^{2}$Institute of Mathematics of the Czech Academy of Sciences, Zitna str., 25, CZ-11567, Prague, Czech Republic,\\
 \indent \,\,\, e-mail: gogatish@math.cas.cz\\
 \indent \,\, \em Manuscript received May 2023}

\begin{abstract}
In this paper, we consider the generalized fractional maximal function and use it to introduce the space of generalized fractional maximal functions and the various cones  of monotone functions generated by generalized fractional maximal functions $M_\Phi f$. We introduced three function classes. We give equivalent descriptions of such cones when the function $\Phi$ belongs to some function classes. The conditions for their mutual covering are given. Then, these cones are used to construct a criterion for embedding the space of generalized fractional maximal functions into the rearrangement invariant spaces (RIS). The optimal RIS for such embedding is also described.


\end{abstract}
\maketitle

\smallskip
\section{Introduction}

In this  paper we consider  the characteristics of function cones exhibiting monotonicity, which are specifically designed for the generalized fractional maximal functions $M_\Phi f$. These functions are defined for $f \in E(\mathbb{R}^n)\cap L_{1}^{loc}(\mathbb{R}^n)$ by

$$(M_{\Phi}f)(x)=\sup_{r>0}{{\Phi(r)}}\displaystyle \int\limits_{B(x,r)} |f(y)|dy,$$
under certain assumptions on the function $\Phi: (0,\infty)\rightarrow (0,\infty)$, where $B(x,r)$ is a ball with the center at the point $x\in\mathbb{R}^n$ and radius $r$, $E=E(\mathbb{R}^{n})$ is the rearrangement-invariant space.

In the case $\Phi(r)=r^{\alpha-n}$, $\alpha\in(0;n)$ we obtain the classical fractional maximal function $M_{\alpha}f$:
$$(M_{\alpha}f)(x)=\sup_{r>0}\frac{1}{r^{n-\alpha}}\displaystyle \int\limits_{B(x,r)} |f(y)|dy.$$

In this paper, we develop some of the results of [4, 5].
We introduce the space of generalized fractional maximal functions constructed by a generalized fractional maximal operator. The base space is a rearrangement-invariant space more general than the Lebesgue space $L_p$.

We introduce four types of cones consisting of non-increasing rearrangements of a generalized fractional-maximal function. The conditions for their mutual covering are determined.

 We employ the axioms of Banach function spaces and rearrangement-invariant spaces as outlined in reference [3], and leverage findings from papers [15-16], which introduce the notions of cone embedding and coverings.

Classical fractional maximal functions and the Riesz potential plays an important role In the theory of functions and in the theory of operators in function spaces, in harmonic analysis, potential theory, and PDE [20, 27, 30] and have important applications to the
Navier-Stokes equations (see [31]) and to the Shr\"odinger equations (see [19]).

These operators can be used in the theory of mathematical modeling and in various issues of optimization problems [1-2, 10-11, 24-25, 28-29, 34].

 It's important to recognize that understanding the  characterization of decreasing rearrangements of potentials and generalized fractional-maximal functions is fundamental in obtaining additional insights into their integral properties. The theory of the maximal functions and the theory of Riesz potentials in classical Lebesgue spaces are described in detail in the books [3], [18], [31-32], [35]. Note that the questions of the boundedness of the classical fractional-maximal operator on general Morrey-type spaces were considered in [8], [9].

The different versions of the generalized fractional-maximal function were previously explored in [21], [23], [26-27]. In those studies, the focus was on investigating the boundedness of the generalized fractional-maximal operators across different function spaces.

 The cones generated by non-increasing rearrangements of the generalized Riesz and Bessel potential were considered by M.L. Goldman [14-16] and by others [6], [7], [17]. In these papers the questions of their mutual covering and the questions of embedding the space of the generalized Riesz and Bessel potentials in rearrangement-invariant spaces were considered. Here we prove that the cones generated by the generalized fractional-maximal function are estimated by the corresponding cones from the generalized Riesz potentials, i.e. the cone of monotone functions generated by generalized fractional functions is covered by the corresponding cone of monotone functions generated by generalized Riesz potentials. The  estimates  for non-increasing rearrangements of a generalized fractional-maximal function are given.

During this estimate, we are obtaining nonlinear so-called supremal operators. We can't use dual operators and we need to develop a new approach for investigating such situations. We are using new results developed for supremal operators to obtain a description of the optimal rearrangements-invariant space for the embedding of the cone of monotone functions generated by generalized fractional functions into rearrangements-invariant spaces.

The article is structured as follows. Section 2 contains some preliminary definitions. In Section 3 we introduce four types of cones consisting of non-increasing rearrangements of a generalized fractional-maximal function. The conditions for their mutual covering are determined. Section 4 explores the embedding of the generalized fractional-maximal function space into rearrangement-invariant spaces. This involves embedding a set of decreasing rearrangements into the corresponding rearrangement-invariant space. We provide a detailed account of the most suitable rearrangement-invariant space for such embedding. Additionally, we establish some notation conventions.
 Throughout this work, we use the letters $C$, $C_1$, $C_2$ for a positive constant, independent of appropriate parameters and not necessarily the same at each occurrence. The notation $f(x)\cong g(x)$ means that there are constants $C_1>0$, $C_2>0$ such that $C_1f(t)\leq g(t)\leq C_2f(t), \;\; t\in \mathbb{R}_+.$

\

\section{Preliminary information}

Consider a space $(S, \Sigma, \mu)$ equipped with a measure. Here, $\Sigma$ represents a $\sigma$-algebra of subsets of the set $S$, over which a non-negative, $\sigma$-finite and $\sigma$-additive measure $\mu$. Let $L_{0} = L_{0}(S, \Sigma, \mu)$ denote the collection of $\mu$-measurable real-valued functions $f: S \rightarrow \mathbb{R}$. Additionally, let $L_{0}^{+}$ be a subset of $L_{0}$ comprising non-negative functions:
 $$L_{0}^{+}=\{f\in L_{0}:f\geq 0\}.$$

 By $L_{0}^+ (0,\infty;\downarrow)$ we denote the set of all non-increasing functions from $L_{0}^+ (0,\infty)$.

 In this work, we will use the concepts of a Banach-functional space (briefly: BFS), introduced by C.Bennett, R.Sharpley [3], as well as the concepts of an ideal space (briefly: IS) considered in the book by S.G.Krein, Yu.I.Petunin and E.M.Semenov [22].

\begin{definition}
{\normalfont{[3]}}
A functional norm, denoted as $\rho: L_{0}^{+}\rightarrow [0, \infty]$, satisfies the following conditions for all $ f,g,f_{n}\in L_{0}^{+},\;\;n \in \mathbb{N}$:

$(P1)\,\rho(f)=0 \Leftrightarrow f=0,$   $\mu$-almost everywhere (briefly:\;\;\;$\mu$-a.e.);

\quad{}\quad{}\, $\rho(\alpha f)=\alpha\rho(f),\alpha\geq 0;     \rho(f+g)\leq \rho(f)+\rho(g)$ (properties of the norm);

$(P2)\, f\leq g,$ $(\mu-$ a.e.) $\Rightarrow \rho(f)\leq\rho(g)$ (monotony of the norm);

$(P3)\, f_{n}\uparrow f  \Rightarrow \rho(f_{n})\rightarrow\rho(f)  (n\rightarrow \infty)$ (the Fatou property);

$(P4)\, 0<\mu(\sigma)<\infty\Rightarrow    \int\limits_{\sigma}f d\mu\leq c_{\sigma}\rho(f), f\in L_{0}^{+}.$ (Local integrability);

$(P5)\, 0<\mu(\sigma)<\infty\Rightarrow   \rho(\chi_{\sigma})<\infty $ (finiteness of the FN for characteristic functions $(\chi_{\sigma})$ of sets of finite measure).

Here $\,f_{n}\uparrow f\,$  means that $f_{n}\leq f_{n+1},  \lim\limits_{n\rightarrow \infty} f_{n}=f\,$ ($\mu$-a.e.).
\end{definition}

\begin{definition}
Let $\rho$ represent a functional norm. A Banach function space (abbreviated as BFS), denoted as $X=X(\rho)$, consists of functions from $L_{0}$ such that $\rho(|f|)<\infty$. For any $f\in X$, we define the norm
$$\|f\|_{X}=\rho(|f|).$$
\end{definition}

Let $\mathbb{R}_+:=(0,+\infty)$ and $L_{0}=L_{0}(\mathbb{R}^{n})$ be the set of all Lebesgue measurable functions $f:\mathbb{R}^{n}\rightarrow \mathbb{C}$. The non-increasing rearrangement $f^{*}$ of a measurable function $f$ on $\mathbb{R}^n$ is defined as
$$ f^{*}(t)=inf\;\{y>0: \;\lambda_{f}(y)\leq t\}, \;t\in\mathbb{R}_{+},$$
where
$$\lambda_{f}(y)=\mu_{n}\left\{x\in \mathbb{R}^{n}: \, |f(x)|>y\right\}, \; y\in\mathbb{R}_{+}$$
is the Lebesgue distribution function, $\mu_{n}$ is Lebesgue measure in $\mathbb{R}^n$ [3].
It is known that $0\leq f^{*}\downarrow;\;\;\; f^{*}(t+0)=f^{*}(t),\;\;\;t\in \mathbb{R}_{+};$ $f^{*}$ is equimeasurable with $f$, i.e.
$$\mu_{1}\left\{t\in R_{+}:\,\,\, f^{*}(t)>y\right\}=\mu_{n}\left\{x\in\mathbb{R}^{n}:\,\,\, |f(x)|>y\right\},$$
here $\mu$ is the Lebesgue measure in $\mathbb{R}^{n}$ or $\mathbb{R}_{+}$, respectively.

The function $f^{**}$ is defined as
$$f^{**}(t)=\frac{1}{t}\int\limits_0^t f^*(\tau)d\tau; \;\; t\in \mathbb{R}_+.$$

Let $f^{\#}: \mathbb{R}^{n}\rightarrow\mathbb{R}^{n}$ represent a symmetrical rearrangement of $f$, which means it is a radially symmetric, non-negative, decreasing, right-continuous function (with respect to $r=|x|$, where $x\in \mathbb{R}^{n}$) that has the same measure as $f$. That is
$$f^{\#}(r)=f^{*}(v_{n}r^{n}); \; f^{*}(t)=f^{\#}\left(\left(\frac{t}{v_{n}}\right)^{\frac{1}{n}}\right),\; \; r, t\in\mathbb{R}_{+}.$$

Here $v_n$ is the volume a unit ball in $\mathbb{R}^{n}$.

\begin{definition}
A functional norm $\rho$ is considered rearrangement-invariant if the inequality $f^{*}\leq g^{*}$ implies $\rho(f)\leq \rho(g)$. A Banach function space $X=X(\rho)$, formed by a rearrangement-invariant functional norm $\rho$, is termed a rearrangement-invariant space, abbreviated as RIS.
\end{definition}

The Lebesgue spaces $L_{p}(\mathbb{R}^{n})$, as well as the Lorentz and Orlich spaces are examples of rearrangement-invariant spaces.

Everywhere in this work, we denote by $E=E(\mathbb{R}^{n})$ the rearrangement-invariant space and by $E^{'}=E^{'}(\mathbb{R}^{n})$ we denote the associated rearrangement invariant space. That is
\begin{equation*}
\|f\|_{E'(\mathbb{R}^{n})}=sup\left\{\left|\int\limits_{R^{n}}f\,gd\mu\right|:\;\;f\in E(\mathbb{R}^{n}),\; \|f\|_{E(\mathbb{R}^{n})}\leq 1\right\}<\infty.
\end{equation*}

 By $\tilde{E}=\tilde{E}(R_{+})$ we denote the Luxemburg representation of $E(\mathbb{R}^{n})$, i.e. rearrangement-invariant space such that
\begin{equation}
\|f\|_{E}=\|f^{*}\|_{\tilde{E}}. \label{1}
\end{equation}

\begin{definition}
A function $f:\mathbb{R}_{+}\rightarrow \mathbb{R}_{+}$ is called quasi-decreasing (quasi-increasing) if there exists a positive constant number $C>1$, such that
$$f(t_{2})<Cf(t_{1}) \;\; \text{if $t_{1}<t_{2}$} \;\;\;\;\;\; (f(t_{1})<Cf(t_{2}) \;\; \text{if $t_{1}<t_{2}$}).$$
\end{definition}

\section{Main results. Equivalent descriptions for the cones of rearrangements of generalized fractional maximal function}

We define the following three function classes.

\begin{definition}
{\normalfont{[5]}}
 Let $n\in\mathbb{N}$ and $R\in (0;\infty]$. We say that a function $\Phi:(0;R)\rightarrow \mathbb{R}_{+}$ belongs to the class $A_{n}(R)$ if:
\\ (1) $\Phi$ is decreases and continuous on $(0;R)$;
\\ (2) the function $\Phi(r)r^{n}$ is quasi-increasing on $(0,R).$
\end{definition}

For example, $\Phi(t)=t^{-\alpha}\in A_{n}(\infty), \;\; 0<\alpha<n.$

\begin{definition}
{\normalfont{[5]}}
 Let $n\in\mathbb{N}$ and let $R\in (0;\infty]$. A function $\Phi:(0;R)\rightarrow \mathbb{R}_{+}$ belongs to the class $B_{n}(R)$ if the following conditions hold:
\\ (1) $\Phi$ is decreases and continuous on $(0;R)$;
\\ (2) There exists a constant $C=C(\Phi)>0$ such that
\begin{equation}
\int\limits^{r}_{0} \Phi{(\rho)}\rho^{n-1}d\rho \leq C\Phi(r)r^{n}, \;\; r\in (0,R). \label{2/2.1}
\end{equation}
\end{definition}

For example,
$$\Phi{(\rho)}=\rho^{\alpha-n}\in B_{n}(\infty) \; (0<\alpha<n);\;\; \Phi{(\rho)}=\ln \frac{eR}{\rho}\in B_{n}(R), R\in \mathbb{R}_{+}.$$
\;\;

For $\Phi\in B_{n}(R)$ the following estimate also holds:
\begin{equation}
\int\limits^{r}_{0} \Phi{(\rho)}\rho^{n-1}d\rho\geq n^{-1}\Phi(r)r^{n},\; r\in(0,R). \label{3/2.2}
\end{equation}

Therefore,
\begin{equation}
\int\limits^{r}_{0} \Phi{(\rho)}\rho^{n-1}d\rho\cong \Phi(r)r^{n}, \; r\in(0,R), \label{4/2.3}
\end{equation}

\begin{equation*}
\Phi\in B_{n}(R)\Rightarrow\{0\leq\Phi\downarrow; \; \Phi(r)r^{n}\uparrow,\; r\in (0,R)\}.
\end{equation*}

For each $\alpha\in[1;\infty)$ there exists  $\beta=\beta(\alpha, c, n)\in[1;\infty)$ ($c$ is a constant from (2)) such that
\begin{equation*}
\Big\{\rho, r\in(0;R); \alpha^{-1}\leq \frac{\rho}{r}\leq\alpha\Big\}\Rightarrow \beta^{-1}\leq \frac{\Phi(\rho)}{\Phi(r)}\leq\beta,
\end{equation*}
\begin{equation*}
(2)\Rightarrow \exists\gamma\in(0;n) \;\; \text{such that} \;\;\Phi(r)r^{\gamma}\uparrow \; \; \text{on} \; (0;R).
\end{equation*}

\begin{definition}

{\normalfont{[4]}} Let $R\in (0;\infty]$. We say that a function $\Phi:(0;R)\rightarrow\mathbb{R}_{+}$ belongs to the class $D(R)$ if for some $C=C(\Phi)>0$
\begin{equation}
\int\limits_{0}^{r^n}\frac{ds}{\Phi(s^{1/n})s}\leq\frac{C}{\Phi(r)}, \;\;\; r\in(0;R). \label{5/2.7}
\end{equation}
\end{definition}

For example,
$$\Phi(t)=t^{-n}ln(1+t)^\alpha\in D(\infty) \; (\alpha>0),$$
$$\Phi(t)=t^{\alpha-n}\in D(\infty) \; (0<\alpha<n).$$

Note that relation (5) is equivalent to the inequality (which can be obtained via a variable transformation):

\begin{equation}
\int\limits_{0}^{r}\frac{dt}{\Phi(t)t}\leq{\frac{C_n}{\Phi(r)}}.  \label{6/2.7'}
\end{equation}
where $C_n=\frac{C}{n}$.

\begin{lemma} \normalfont{[4]}. Let $n\in\mathbb{N}$, $R\in(0,\infty]$. Then  $B_{n}(R)\varsubsetneqq A_n(R)$.
\end{lemma}

The function $\Phi(t)=t^{-n}\ln(1+t)^\alpha$, with $\alpha >0$, belongs to $A_n(R)$ and $\Phi\not \in B_n(R) (see [4]).$

We define the space of generalized fractional maximal functions $M_{E}^{\Phi}\equiv M_{E}^{\Phi}(\mathbb{R}^{n})$ as the collection of functions $u$ such that there exists a function $f\in E(\mathbb{R}^{n})$ satisfying
$$u(x)=(M_{\Phi}f)(x),$$
$$\|u\|_{M_{E}^{\Phi}}=inf\{\|f\|_E:f\in E(\mathbb{R}^n), \;\; M_\Phi f=u\}.$$

Note that the generalized Riesz potential was considered in [14-16], [6], [7], [17] using the convolution operator:
$$(I_G)f(x)=(G\ast f)(x)=2\pi^{-n/2}\displaystyle \int\limits_{\mathbb{R}^{n}} G(x-y)f(y)dy,$$
where the kernel $G(x)$ satisfies the conditions:
\begin{equation}
G(x)\cong \Phi(|\rho|), \;\; \rho=|x|\in\mathbb{R}_+, \Phi\in B_n(\infty). \label{7/2.13}
\end{equation}

The classical Riesz potential's kernel takes the shape of
$$G(x)=|x|^{\alpha-n}, \;  \alpha\in(0,n).$$

The following lemma shows that the generalized fractional-maximal function $M_{\Phi}f(x)$ is estimated by the generalized Riesz potential $I_{\Phi}f(x)$. This result is somehow in the proof of Theorem 2 in the paper [4]. As it is not formulated in such a way and it will used many times, we present it as a separate Lemma.

\begin{lemma} {\normalfont{[4]}}. Let $\Phi\in B_{n}(0,\infty)$ and $G(x)= \Phi(|x|), \;\; x\in \mathbb{R}^{n}.$ Then
\begin{equation*}
(M_{\Phi}f)(x)\leq (I_G|f|)(x), \;\; x\in \mathbb{R}^{n}
\end{equation*}
for all $f\in E(\mathbb{R}^n)$.
\end{lemma}

The following three theorems were proved in [4]. We will need them to prove of Theorems 3.5; 3.6; 4.1 and Theorem 4.2.

\begin{theorem} {\normalfont{[4]}} Let $\Phi\in A_{n}(\infty)$. Then there exists a positive constant $C$, depending only on $\Phi$ and $n$, such that
$$(M_{\Phi}f)^{*}(t)\leq C \sup_{t<s<\infty}{s{\Phi}(s^{1/n})} f^{**}(s), \;\;\; t\in(0,\infty),$$
for every $f\in L^1_{loc} (\mathbb{R}^n)$. This inequality is sharp in the sense that for every $\varphi\in L^+_0 (0,\infty;\downarrow)$ there exists a function $f\in L^+(\mathbb{R}^n)$ such that $f^*=\varphi$ a.e. on $(0,\infty)$ and
\begin{equation*}
(M_{\Phi}f)^{*}(t)\geq C_1 \sup_{t<s<\infty}{s{\Phi}(s^{1/n})} f^{**}(s),\;\;\; t\in(0,\infty),
\end{equation*}
where $C_1$ is a positive constant which depends only on $\Phi$ and $n$.
\end{theorem}
For the case of classical fractional maximal function, such theorem was proved in [12].
\begin{theorem} {\normalfont{[4]}} Let $\Phi\in B_{n}(\infty)$. Then there exists a positive constant $C$, depending only on $\Phi$ and $n,$ such that
\begin{equation*}
(M_{\Phi}f)^{**}(t)\leq C \sup_{t<s<\infty}{s{\Phi}(s^{1/n})} f^{**}(s), \;\;\; t\in(0,\infty),
\end{equation*}
for every $f\in L^1_{loc} (\mathbb{R}^n)$.
\end{theorem}

\begin{theorem} {\normalfont{[4]}} Let $\Phi\in B_n(\infty)\cap D(\infty)$, then for every $f\in L^1_{loc}(\mathbb{R}^n)$ there exists a positive constant $C$, depending only on $\Phi$ and $n$, such that
\begin{equation*}
(M_{\Phi}f)^{*}(t)\leq C \Bigg(t{\Phi}(t^{1/n})f^{**}(t)+ \sup_{t<\tau<\infty} {\tau{\Phi}(\tau^{1/n})} f^{*}(\tau)\Bigg), \;\;\; t\in(0,\infty).
\end{equation*}
\end{theorem}

\begin{definition}

Let $\Im_T={K(T)}$ for $T\in(0,\infty]$ be defined as the collection of cones containing measurable non-negative functions defined on $(0,T)$, together with positive homogeneous functionals $\rho_{K(T)}:K(T)\rightarrow [0,\infty)$ satisfying the following properties:

$(1) \;\; h\in K(T), \;\; \alpha\geq 0\Rightarrow \alpha h\in K(T), \;\;\;\; \rho_{K(T)}(\alpha h)=\alpha\rho_{K(T)}(h);$

$(2) \;\; \rho_{K(T)}(h)= 0 \Rightarrow h=0$ almost everywhere on $(0,T)$.
\end{definition}

\begin{definition} {\normalfont{[15]}} Let $K(T), M(T)\in\Im_{T}$. The cone $K(T)$ {\itshape covers the cone}  $M(T)$ (notation: $M(T)\prec K(T)$) if there exist $C_{0}=C_{0}(T)\in R_{+}$, and $C_{1}=C_{1}(T)\in [0,\infty)$ with $C_{1}(\infty)=0$ such that for each $h_{1}\in M(T)$ there is  $h_{2}\in K(T)$ satisfying
$$\rho_{K(T)}(h_{2})\leq C_{0}\rho_{M(T)}(h_{1}), \;\; h_{1}(t)\leq h_{2}(t)+C_{1}\rho_{M(T)}(h_{1}), \;\; t\in(0,T).$$
\end{definition}

The equivalence of the cones means mutual covering:
$$M(T)\approx K(T)\Leftrightarrow M(T)\prec K(T)\prec M(T).$$

Let $E$ be a rearrangement-invariant space. We examine the following four cones generated by the non-increasing rearrangements of the generalized fractional maximal function, each equipped with homogeneous functionals:
\begin{equation}
K_1\equiv K_{E}^{\Phi}:=\{h\in L^{+}(\mathbb{R_+}):h(t)=u^{*}(t), \; t\in \mathbb{R_+}, \; u\in M_{E}^{\Phi}\}, \label{8/3.3}
\end{equation}
$$\rho_{K_1}(h)=\inf\{\|u\|_{M_{E}^{\Phi}}: u\in M_{E}^{\Phi}; \; u^{*}(t)=h(t), \; t\in \mathbb{R_+}\};$$

\begin{equation}
K_2\equiv\widehat{K_E^\Phi}:=\{h:h(t)=u^{**}(t), \; t\in\mathbb{R_+}, \; u\in M_{E}^{\Phi}\}, \label{9/3.5}
\end{equation}
$$\rho_{K_2}(h)=\inf\{\|u\|_{M_{E}^{\Phi}}:u\in M_{E}^{\Phi}; \; u^{**}(t)=h(t), \; t\in \mathbb{R_+}\};$$

\begin{equation}
K_3\equiv\widetilde{K_E^\Phi}:=\big\{h\in L^{+}(\mathbb{R_+}):h(t)= \sup_{t<\tau<\infty}\tau\Phi(\tau^{1/n})u^{**}(\tau),\; t\in\mathbb{R_+},\; u\in E(\mathbb{R}^n)\big\},  \label{10/3.7}
\end{equation}

\begin{equation}
\rho_{K_3}(h)=inf\big\{\|u\|_{E(\mathbb{R}^n)}, u\in E(\mathbb{R}^n): \sup_{t<\tau<\infty}\tau\Phi(\tau^{1/n})u^{**}(\tau)=h(t),\; t\in\mathbb{R_+}\big\}; \label{11/3.8}
\end{equation}

\begin{equation}
K_4\equiv\widetilde{\widetilde{K^\Phi_E}}:=\big\{h\in L^{+}(\mathbb{R_+}):h(t)= t\Phi(t^{1/n})u^{**}(t) +\sup_{{t<\tau<\infty}} \tau\Phi(\tau^{1/n})u^{*}(\tau), t\in\mathbb{R_+},\; u\in E(\mathbb{R}^n)\big\}, \label{12/3.9}
\end{equation}

\begin{equation}
\rho_{K_4}(h)=inf\big\{\|u\|_{E(\mathbb{R}^n)}, u\in E(\mathbb{R}^n): t\Phi(t^{1/n})u^{**}(t)+\sup_{t<\tau<\infty}\tau\Phi(\tau^{1/n})u^{*}(\tau)=h(t),\; t\in\mathbb{R_+}\big\}. \label{13/3.10}
\end{equation}

Aware that Goldman M.L. [14-16] (as well as references [6]-[7], [17]) delve into the examination of cones generated by generalized potentials. They investigate the realm of generalized Riesz potentials denoted as $H_{E}^{G}(\mathbb{R}^n)$ within $n$-dimensional Euclidean space:
$$H_{E}^{G}(\mathbb{R}^n)=\{u=G*f:f\in E(\mathbb{R}^n)\},$$
$$\|u\|_{H_{E}^{G}}=\inf\{\|f\|_{E}: f\in E(\mathbb{R}^n); G*f=u\},$$
where $E(\mathbb{R}^n)$ is an rearrangement invariant space. In that work, the following cones of rearrangements for generalized Riesz potential equipped with a positive homogeneous functional were considered:
$$M\equiv M_{E,G}(\mathbb{R_+})=\{h(t)=u^{*}(t), t\in (0;\infty), u\in H_{E}^{G}\},$$
$$\rho_{M}(h)=\inf\{\|u\|_{H_{E}^{G}}: u\in H_{E}^{G}; u^{*}(t)=h(t), t\in (0;\infty)\};$$
$$\widetilde{M}\equiv\widetilde{M}_{E,G}(\mathbb{R_+})=\{h(t)=u^{**}(t), t\in (0;\infty): u\in H_{E}^{G}\},$$
$$\rho_{\widetilde{M}}(h)=\inf\{\|u\|_{H_{E}^{G}}:u\in H_{E}^{G}; u^{**}(t)=h(t), t\in (0;\infty)\}.$$

For $\tau, t\in R_+$ denote
$$\varphi(\tau)=\Phi\Big(\big(\frac{\tau}{v_n}\big)^{\frac{1}{n}}\Big)\in B_1(\infty), f_{\Phi}(t;\tau)=\min\{\varphi(t),\varphi(\tau)\}.$$

In [15] was also considered the cone
$$K=K_{E,\Phi}(\infty)=\bigg\{h(g;t)=h(t)=\int\limits_{0}^{\infty} f_{\Phi}(t;\tau)g(\tau)d\tau, t\in\mathbb{R_+}: g\in{\widetilde{E}_{0}(\mathbb{R_+})}\bigg\}  $$
equipped with a functional

$$\rho_{K}(h)=inf\big\{\|g\|_{\widetilde{E}(\mathbb{R_+})}:h(g;t)=h(t),t\in\mathbb{R_+}\big\},$$
where
\begin{equation}
\tilde{E}_{0}(0,\infty)=\big\{g\in{\widetilde{E}(\mathbb{R_+})}: 0\leq g\downarrow; g(t+0)=g(t), t\in (\mathbb{R_+})\big\}.  \label{14/3.13}
\end{equation}

In [15] was given the equivalent characterization for these cones of decreasing rearrangements of potentials.

The following theorem proves that the cone $M$ generated by the generalized Riesz potential covers the cone generated by the generalized maximal function.

\begin{theorem} {\normalfont{[5]}} Let $\Phi\in B_n(\infty)$ and kernel $G(x)$ satisfies the condition (7). Then the cone generated by the generalized potential covers the cone generated by the generalized maximal function, i.e. $K_1\prec M$.
\end{theorem}

The following theorem gives an equivalent description of the cones $K_1, K_2, K_3$.

\begin{theorem} Let $\Phi\in B_{n}(\infty)$. Then cones (8), (9), (10) are equivalent:
$$K_1\approx K_2\approx K_3.$$
\end{theorem}
\emph{Proof.}

1) It is easy to check that  $K_1\prec K_2$ (see [5]).

2) Prove that $K_2\prec K_3$. Let $h_3\in K_2$, then there is $u_2\in {M}_{E}^{\Phi}$ such that

$$u_2^{**}=h_{3}, \;\;\; \|u\|_{M_{E}^{\Phi}}\leq 2\rho_{K_2}(h_3).$$

For $u_2\in M_{E}^{\Phi}$ we find $f_2\in E(R^n)$ satisfying the condition

$$u_2(x)=(M_{\Phi}f_2)(x)=\sup_{r>0}\Phi(r)\int\limits_{B(x,r)} |f(\xi)|d\xi,$$

$$\|f_2\|_{E(\mathbb{R}^n)}\leq 2\|u_2\|_{M_{E}^{\Phi}}.$$

Hence
\begin{equation}
h_3(t)=(M_{\Phi}f_2)^{**}(t), \;\; \|f_2\|_{E(\mathbb{R}^n)}=\|f_2^*\|_{\widetilde{E}(R_+)}\leq 4 \rho_{K_2}(h_3). \label{15/3.36}
\end{equation}

According to Theorem 3.2,
\begin{equation}
(M_{\Phi}f_2)^{**}(t)\leq C \sup_{t<\tau<\infty}{\tau{\Phi}(\tau^{1/n})}f_2^{**}(\tau).\label{16/3.37}
\end{equation}

We denote
$$h_4=C\sup_{t<\tau<\infty} {\tau{\Phi}(\tau^{1/n})}f_2^{**}(\tau).$$

Hence, by (15), (16) we have $h_3(t)<h_4(t)$ and

\begin{equation*}
\rho_{K_3}(h_4)=inf\big\{\|u\|_{E(\mathbb{R}^n)}: \sup_{t<\tau<\infty}\tau\Phi(\tau^{1/n})u^{**}(\tau)=h_4(t)\big\}\leq\|f_2\|_{E(\mathbb{R}^n)}= \|f_2^*\|_{\widetilde{E}(R^+)}\leq 4\rho_{K_2}(h_3).
\end{equation*}

Covering $K_2\prec K_3$ is proved.

3) Now we prove that the embedding $K_3\prec K_1$ takes place when $\Phi\in B_{n}(\infty)$. Let $h_5\in K_3$. Then there is $u_3\in E(\mathbb{R}^n)$ such that
$$h_5(t)=\sup_{t<\tau<\infty} {\tau{\Phi}(\tau^{1/n})}u_3^{**}(\tau), t\in\mathbb{R}_+,$$
$$\|u_3\|_{E(\mathbb{R}^n)}\leq2\rho_{K_3}(h_5).$$

According to Theorem 3.1 for $u_3\in L_0^+(0,\infty;\downarrow)$ there exist the function $f_3$ on $\mathbb{R}^n$ such that $f_3^*=u_3^*$ a.e. on $(0,\infty)$ and
$$(M_\Phi f_3)^*(t)\geq C\sup_{t<\tau<\infty} {\tau{\Phi}(\tau^{1/n})}f_3^{**}(\tau), t\in\mathbb{R}_+.$$

We denote
$$h_6(t)=(M_\Phi f_3)^*(t)\in K_1.$$

Therefore,
$$h_5(t)\leq \frac{1}{C}h_6(t)$$
$$\rho_{K_1}(h_6)\leq\|u_3\|_{M_E^\Phi}\leq\|f_3^*\|_{\widetilde{E}(R^+)}=\|u_3\|_{E(\mathbb{R}^n)}\leq2\rho_{K_3}(h_5).$$

So we proved $K_3\prec K_1$.

Theorem 3.5 is proved.

\begin{theorem}
 1) Let $\Phi\in B_{n}(\infty)$. Then for cones (10), (12)
 $$K_4\prec K_3$$
\\
 2)Let $\Phi\in D(\infty)$. Then for cones (10), (12)
 $$K_3\prec K_4.$$
\\
 3) Let $\Phi\in B_n(\infty)\cap D(\infty)$. Then for cones (10), (12)
 $$K_3\approx K_4.$$
\end{theorem}

\emph{Proof.}
 1) Let $h_1\in K_4$. Then there is a function $u_1\in E(\mathbb{R}^n)$ such that
\begin{equation}
h_1(t)=t{\Phi}(t^{1/n})u_1^{**}(t)+\sup_{t<\tau<\infty}{\tau{\Phi}(\tau^{1/n})}u_1^*(\tau)
\label{17/3.20}
\end{equation}
and (see (1), (13))
\begin{equation}
\|u_1\|_{E(\mathbb{R}^n)}=\|u_1^*\|_{\widetilde{E}_0(\mathbb{R}_+)}\leq2\rho_{K_4}(h_1). \label{18/3.21}
\end{equation}
By inequality
$$u_1^*(\tau)\leq u_1^{**}(\tau), \; \tau\in\mathbb{R}_+,$$

from (17) we have
$$h_1(t)\leq t{\Phi}(t^{1/n})u_1^{**}(t)+ \sup_{t<\tau<\infty}{\tau{\Phi}(\tau^{1/n})}u_1^{**}(\tau)= \sup_{t<\tau<\infty}{\tau{\Phi}(\tau^{1/n})}u_1^{**}(\tau).$$

We denote
$$h_2(t)=\sup_{t<\tau<\infty}{\tau{\Phi}(\tau^{1/n})}u_1^{**}(\tau)\in K_3.$$

Then
$$h_1(t)\leq h_2(t), \; t\in\mathbb{R}_+.$$

Moreover, (see (18), (11))
$$\rho_{K_3}(h_2)\leq\|u_1\|_{E(\mathbb{R}^n)}\leq2\rho_{K_4}(h_1).$$

\emph{Item 1 is proved.}

2) Let $h_3\in K_3$. Then there is a function $u_2\in E(\mathbb{R}^n)$ such that
$$h_3(t)=\sup_{t\leq\tau<\infty}{\tau{\Phi}(\tau^{1/n})}u_2^{**}(\tau)$$
and (see (1), (13))
\begin{equation}
\|u_2\|_{E(\mathbb{R}^n)}\leq2\rho_{K_3}(h_3). \label{19/3.25}
\end{equation}

According to Theorem 3.3
$$h_3(t)=\sup_{t\leq\tau<\infty}{\tau{\Phi}(\tau^{1/n})}u_2^{**}(\tau)\leq C_1 \bigg(t{\Phi}(t^{1/n})u_2^{**}(t)+\sup_{t<\tau<\infty}{\tau{\Phi}(\tau^{1/n})}u_2^*(\tau)\bigg).$$

Let's put
$$h_4(t)=C_1 \bigg(t{\Phi}(t^{1/n})u_2^{**}(t)+ \sup_{t<\tau<\infty}{\tau{\Phi}(\tau^{1/n})}u_2^*(\tau)\bigg)\in K_4.$$

So
$$h_3(t)\leq h_4(t), t\in \mathbb{R}_+.$$

Moreover (see (13), (19))
$$\rho_{K_4}(h_4)\leq\|u_2\|_{E(\mathbb{R}^n)}\leq2\rho_{K_3}(h_3).$$
\emph{Item 2 is proved.}

3) Item 3 follows from items 1 and 2.

Theorem 3.6 is proved.

\begin{corollary}
Let $\Phi\in B_n(\infty)\cap D(\infty)$. Then
$$K_1\approx K_2\approx K_3\approx K_4.$$
\end{corollary}
The Corollary 3.1 follows from Lemma 3.1, Theorem 3.5, and Theorem 3.6.

\section{On embedding of spaces of generalized fractional maximal functions in rearrangement invariant spaces}

Let $T\in(0,\infty]$, $\Im_{T}=\{M(T)\}$ the set of cones consisting measurable function on $(0,T)$ equipped with a positive homogeneous functionals $\rho_{M(T)}$ (see Definition 3.5). Let $X(\mathbb{R}^n)$ is the rearrangement invariant space, $\widetilde{X}(R^+)$ is their Luxemburg representation (1).

\begin{definition} {\normalfont{[14]}}
Let $M(T)\in\Im_{T}$. The embedding $M(T)\mapsto\widetilde{X}(0,T)$ means that $M(T)\subset\widetilde{X}(0,T)$ and there exists a constant $C_{M(T)}(T)\in R_{+}$ such that
$$\|h\|_{\widetilde{X}(0,T)}\leq C_{M(T)}\rho_{M(T)}(h), \;\; h\in M(T).$$

Allow $M(T)$ to belong to $\Im_{T}$. The mapping $M(T)\mapsto\widetilde{X}(0,T)$ implies that $M(T)$ is contained within $\widetilde{X}(0,T)$, and there exists a positive constant $C_{M(T)}(T)\in \mathbb{R}_{+}$ such that for any $h\in M(T)$, we have

\end{definition}

\begin{lemma} {\normalfont{(see [14])}}
Let $K(T), M(T)\in\Im_{T}$. For every Banach function space $\widetilde{X}(0,T)$ if $M(T)\prec K(T)$ then $K(T)\mapsto\widetilde{X}(0,T)\Rightarrow M(T)\mapsto\widetilde{X}(0,T)$.
\end{lemma}

\begin{corollary}
Let $K(T), M(T)\in\Im_{T}$. If $M(T)\approx K(T)$ then $K(T)\mapsto\widetilde{X}(0,T)\Leftrightarrow M(T)\mapsto\widetilde{X}(0,T)$.
\end{corollary}

Let us formulate criteria for embedding of the space of generalized fractional maximal functions in rearrangement invariant spaces $X(\mathbb{R}^n)$:
\begin{equation}
M_E^\Phi(\mathbb{R}^n)\hookrightarrow X(\mathbb{R}^n).  \label{20/4.3}
\end{equation}

\begin{definition}
Let $E$ - RIS. By an optimal RIS for embedding (20) we mean an RIS $X_0=X_0(R^n)$ such that $M_E^\Phi(\mathbb{R}^n)\mapsto X_0(\mathbb{R}^n)$ and if (20) is valid for another RIS $X$, then $X_0\subset X$. Such an optimal RIS is called a rearrangement-invariant envelope of the space of fractional-maximal functions.

Let $E$ be a rearrangement-invariant space (RIS). An optimal RIS for embedding (20) is defined as an RIS $X_0 = X_0(\mathbb{R}^n)$ such that the embedding operator $M_E^\Phi(\mathbb{R}^n) \rightarrow X_0(\mathbb{R}^n)$ and if (20) holds for another RIS $X$, then $X_0$ is a subset of $X$. This optimal RIS is referred to as a rearrangement-invariant envelope of the space of fractional-maximal functions.
\end{definition}

\begin{theorem} {\normalfont{[5]}} Let $\Phi\in B_n(\infty)$. Then the embedding (20) is equivalence to the following embedding
\begin{equation*}
K_1\mapsto \widetilde{X}(\mathbb{R}_+).
\end{equation*}
\end{theorem}

The next corollary follows from Theorem 3.5, Theorem 4.1, and Corollary 4.1.

\begin{corollary}

1) Let  $\Phi\in B_n(\infty)$. Then
\begin{equation*}
M_E^\Phi(\mathbb{R}^n)\subset X(\mathbb{R}^n) \Leftrightarrow K_i\mapsto \widetilde{X}(\mathbb{R}_+), \;\; i=1;2;3.
\end{equation*}

2) Let $\Phi\in B_n(\infty)\cap D(\infty)$. Then
\begin{equation*}
M_E^\Phi(\mathbb{R}^n)\subset X(\mathbb{R}^n) \Leftrightarrow K_i\mapsto \widetilde{X}(\mathbb{R}_+), \;\; i=1;2;3;4.
\end{equation*}
\end{corollary}

We denote
$$(T_\Phi f)(x)=\sup_{x<t<\infty}\Phi(t^{1/n})tf^{**}(t),$$
operator $T_\Phi: \widetilde{E}_0(\mathbb{R}_+)\rightarrow {X_0}(\mathbb{R}_+)$, where $\widetilde{E}_0(\mathbb{R}_+)$ defined by (14).

\begin{theorem}
Let $\Phi\in B_n(\infty)$. Then the embedding
\begin{equation*}
K_3\mapsto \widetilde{X}(\mathbb{R}_+)
\end{equation*}
is equivalent to the boundedness of the operator $T_\Phi$ from $\widetilde{E}(\mathbb{R}_+)$ to $\widetilde{X}(\mathbb{R}_+)$.
\end{theorem}

\emph{Proof.} Indeed, this embedding  means that

$$h(t)=T_\Phi f(t)=\sup_{t<s<\infty}\Phi(s^{1/n})sf^{**}(t)\in\widetilde{X}(\mathbb{R}_+), \; f^*\in\widetilde{E}(\mathbb{R}_+),$$
$$\|h\|_{\widetilde{X}(\mathbb{R}_+)}\leq C_0\rho_{K_3}(h)\leq C_0\|f^*\|_{\widetilde{E}(\mathbb{R}_+)},$$
$$\|T_\Phi f\|_{\widetilde{X}(\mathbb{R}_+)}\leq C_0\|f^*\|_{\widetilde{E}(\mathbb{R}_+)},$$
which is equivalent to the boundedness of the operator $T_\Phi$.

By Theorem 3.2 and Theorem 4.1 we have
$$K_3(\mathbb{R}_+)\mapsto \widetilde{X}(\mathbb{R}_+)\Leftrightarrow K_1(\mathbb{R}_+)\mapsto \widetilde{X}(\mathbb{R}_+)\Leftrightarrow M_E^\Phi(\mathbb{R}^n)\hookrightarrow X(\mathbb{R}^n).$$

Next, we need the following result.

\begin{theorem} {\normalfont{[30]}} Suppose $f, w\in L_0^+$. Then
$$\sup_{g:\int_{0}^{x}g\leq\int_{0}^{x}f}\int_{0}^{\infty}g(x)w(x)dx=\int_{0}^{\infty}f(x)\bigg(\sup_{t\in[x,\infty)}w(t)\bigg)dx.$$
\end{theorem}

\begin{theorem}
Let $\Phi\in B_n(\infty)$. The optimal RIS $X_0=X_0(\mathbb{R}^n)$ for embedding
$$M_E^\Phi(\mathbb{R}^n)\hookrightarrow X(\mathbb{R}^n)$$
possesses an equivalent norm, as expressed in the Luxemburg representation (1):
\begin{equation*}
\|f\|_{\widetilde{X}_0(0,\infty)}=\sup_{g^*}\bigg\{{\int\limits_{0}^{\infty}f^{*}(\tau)g^{*}(\tau)d\tau: g\in L_0(0,\infty), \sup_{\int\limits_{0}^{t}h(s)ds\leq\int\limits_{0}^{t}g^*(s)ds} \big\|\int\limits_{t}^{\infty}{\Phi}(s^{1/n})}sh(s)ds\big\|_{E'}\leq1\bigg\}.
\end{equation*}

\end{theorem}

\emph{Proof.}
According to Theorem 4.1 and Theorem 4.3, the embedding (20) is equivalent to the boundedness of the operator $T_\Phi$ from $\widetilde{E}_0(\mathbb{R}_+)$ to ${\widetilde{X}}(\mathbb{R}_+)$, i.e. let us prove the inequality:

$$\|T_\Phi f\|_{X_0}\leq C\|f\|_E.$$

Using the associated space $X_0'$ for the space $X_0$ and using the Theorem 4.3, we have
\begin{eqnarray*}
\sup_{f}\frac{\|T_\Phi f\|_{X_0}}{\|f\|_E} &=& \sup_{f}\sup_{g}\frac{\int\limits_{0}^{\infty}(T_\Phi f)(\tau)g^*(\tau)d\tau}{ \|f\|_E\|g\|_{X'_0}} \\
&=&\sup_{f}\sup_{g}\frac{\int\limits_{0}^{\infty}g^*(\tau)\sup\limits_{\tau<t<\infty}\Phi(t^{1/n})tf^{**}(t)d\tau}{\|f\|_E\|g\|_{X'_0}} \\
&=& \sup_{g} \frac{1}{\|g\|_{X'_0}} \sup_{f} \frac{\int\limits_{0}^{\infty}g^*(\tau){\sup_{\tau<t<\infty}} \Phi(t^{1/n})tf^{**}(t)d\tau}{\|f\|_E} \\
&=&\sup_{g}\frac{1}{\|g\|}_{X'_0}\sup_{f} \sup\limits_{\int\limits_{0}^{t}h(s)ds\leq\int\limits_{0}^{t}g^*(s)ds} \frac{\int_{0}^{\infty}f^*(t)\int_{t}^{\infty}\Phi(s^{1/n})sh}{\|f\|_E}\\
&=&\sup_{g}\frac{1}{\|g\|}_{X'_0} \sup_{\int\limits_{0}^{t}h(s)ds\leq\int\limits_{0}^{t}g^*(s)ds} \left\|\int\limits_{t}^{\infty}{\Phi}(s^{1/n})sh(s)ds\right\|_{E'}.
\end{eqnarray*}

Using the definition of the associate norm we obtain the proof of the Theorem 4.4.

Note that the proof of Theorem 4.4 contains the formula for associate norms of optimal
embedding (20). We have the following Corollary.

\textbf{Corollary 4.1}. Let $\Phi\in B_n(\infty)$. Let $X_0=X_0(\mathbb{R}^n)$ be an optimal RIS space for the embedding (20). Then we have the following formula for the associate space $X_0'$ :
\begin{equation*}
\|g\|_{X'_0}= \sup_{f} \frac{\int\limits_{0}^{\infty}g^*(\tau)\sup\limits_{\tau<t<\infty} \Phi(t^{1/n})tf^{**}(t)d\tau}{\|f\|_E}.
\end{equation*}

\

\section{Conclusion}

In this paper, we considered the space of generalized fractional maximal functions and
investigated the various cones generated by the nonincreasing rearrangement of generalized
fractional maximal functions. The equivalent descriptions of such cones and the conditions for their mutual covering are given. Then these cones are used to construct a criterion for embedding the space of generalized fractional maximal functions into the rearrangement invariant spaces (RIS). The optimal RIS for such embedding is also described.

\

\section{Acknowledgement}

\noindent{The research of A.N.~Abek, N.A.~Bokayev, A.~Gogatishvili was supported by the grant Ministry of Education and Science of the Republic of Kazakhstan (project no: AP14869887).}

\noindent{The research of A. Gogatishvili was  supported by the Czech Academy of Sciences RVO: 67985840, by Czech Science Foundation (GA\v CR), grant no: 23-04720S, and by Shota Rustaveli National Science Foundation (SRNSF), grant no: FR21-12353.}

\end{document}